\documentclass[10pt]{amsart}
\usepackage[left=2.7cm,right=2.7cm,top=3cm,bottom=3cm]{geometry}

\usepackage[colorlinks=true,urlcolor=blue,
citecolor=red,linkcolor=blue,linktocpage,pdfpagelabels,
bookmarksnumbered,bookmarksopen]{hyperref}

\usepackage{amsmath,amssymb,cite}
\usepackage{amstext}
\usepackage{enumerate}
\usepackage{color}

\newtheorem{innercustomthm}{Lemma}
\newenvironment{customthm}[1]
  {\renewcommand\theinnercustomthm{#1}\innercustomthm}
  {\endinnercustomthm}

\def\r{\mathbb{R}}

\def\rn{\mathbb{R}^N}

\def\ch{\mathcal{H}}

\def\eps{\varepsilon}

\def\irn{\int_{\rn}}

\def\wt{\widetilde}

\def\d12{\mathcal{D}^{1,2}}

\title[Logarithmic NLS with periodic potential]{Multiple solutions to logarithmic Schr\"odinger \\ equations
	with periodic potential}

\author[M.\ Squassina]{Marco Squassina}
\address{Dipartimento di Matematica e Fisica
\newline\indent
Universit\`a Cattolica del Sacro Cuore
\newline\indent
Via dei Musei 41, 25121 Brescia, Italy }
\email{marco.squassina@unicatt.it}

\author[A.\ Szulkin]{Andrzej Szulkin}
\address{Department of Mathematics 
\newline\indent
Stockholm University
\newline\indent
106 91 Stockholm, Sweden}
\email{andrzejs@math.su.se}

\begin{document}

\begin{abstract}
In this note we provide a corrigendum to some results in \href{http://arxiv.org/abs/1404.4790}{arXiv:1404.4790}.
\end{abstract}

\maketitle

In our paper \cite{sqsz} there are some technical errors which do not affect the validity of the results but require modification of the arguments. We are grateful to Kazunaga Tanaka, Chengxiang Zhang and Qingye Zhang for kindly pointing out these errors.

We shall use the same notation as in \cite{sqsz}. Lemma 2.4 is not correct as stated and requires some changes. For $u\in D(J)$ and $z \in C_0^\infty(\rn)$, let 
\[
\langle J'(u),z\rangle :=\langle\Phi'(u),z\rangle+\irn Q(x)F_1'(u)z\,dx.
\]
$J'(u)$ is a densely defined linear operator and we may define
\[
\|J'(u)\| := \sup\{\langle J'(u),z\rangle: z\in C_0^\infty(\rn),\ \|z\|\le 1\}.
\]
If $\|J'(u)\|$ is finite, then $J'(u)$ may be extended to a bounded operator in $E$ and is thus an element of $E'$. 
\begin{customthm}{2.4} \label{lem2.4}
Let $u\in D(J)$ and $J'(u)\in E'$. Then $w\in \partial J(u)$, i.e.
\begin{equation} \tag{$*$} \label{subdiff}
\langle \Phi'(u),v-u\rangle+\Psi(v)-\Psi(u) \ge \langle w, v-u\rangle \quad \text{for all $v\in E$}
\end{equation} 
 if and only if $w=J'(u)$. If $(u_n)$ is a Palais-Smale sequence, then $J'(u_n)\in E'$ and $J'(u_n)\to 0$.
\end{customthm} 

\begin{proof}
Insert $v=u+tz$, where $t>0$ and $z\in C_0^\infty(\rn)$, into \eqref{subdiff}. This gives
\[
\langle\Phi'(u),z\rangle + \irn Q(x)\,\frac{F_1(u+tz)-F_1(u)}t\,dx \ge \langle w,z\rangle.
\]
Since the integrand is 0 for $x\not\in\text{supp\,}z$, we can pass to the limit as $t\to 0$ and we obtain
\[
\langle\Phi'(u),z\rangle+\irn Q(x)F_1'(u)z\,dx \ge \langle w,z\rangle
\]
(recall that by Lemma 2.2, $\Psi\in C^1(H^1(\Omega), \r)$ if $\Omega$ is bounded). Since this also holds for $-z$,
\[
\langle J'(u),z\rangle = \langle\Phi'(u),z\rangle+\irn Q(x)F_1'(u)z\,dx = \langle w,z\rangle \quad \text{for all $z\in C_0^\infty(\rn)$}.
\]
By density of $C_0^\infty(\rn)$ in $E$, $w$ is unique and $J'(u)=w$.

If $(u_n)$ is a Palais-Smale sequence, then by definition
\[
\langle \Phi'(u_n), v-u_n\rangle  + \Psi(v)-\Psi(u_n) \ge -\eps_n\|v-u_n\| \quad \text{for all $v\in E$},
\]
where $\eps_n\to 0^+$. Letting $v=u_n+tz$, where $t>0$ and $z\in C_0^\infty(\rn)$, we have
\[
\langle \Phi'(u_n), z\rangle  + \irn Q(x)\,\frac{F_1(u_n+tz)-F_1(u_n)}t\,dx \ge -\eps_n\|z\|,
\]
so passing to the limit as $t\to 0$, we obtain
\[
\langle J'(u_n),z\rangle \ge -\eps_n\|z\|, \quad z\in C_0^\infty(\rn).
\]
Hence $J'(u_n)\in E'$ and $\|J'(u_n)\|\le\eps_n\to 0$. 
\end{proof}

%Note that in view of the new Lemma 2.4, Definition 2.5 and Corollary 2.6 become superfluous. 

The proof of Lemma 2.7 in \cite{sqsz} contains a gap which may not be easy to remove. Instead we prove a different (and in fact simpler!) version of this lemma. We employ an idea from \cite{dpks}. For $d>0$ and $u\in J^d$, let 
\[
\mu_d(u) := \inf_{a\in J^d}\{\|J'(a)\|+\|u-a\|\}.
\]
It is easy to see that $\mu_d(u)<\infty$. Since $\mu_d(u)\le \|J'(u)\|$, it is clear that $\mu_d(u_n)\to 0$ if $J'(u_n)\to 0$. On the other hand, if $\mu_d(u_n)\to 0$, then $J'(a_n)\to 0$ for some $(a_n)$ such that $u_n-a_n\to 0$. The mapping $\mu_d$ is Lipschitz continuous. Indeed, for $u,v,a\in J^d$,
\[
\mu_d(u) \le \|J'(a)\| + \|u-a\| \le \|J'(a)\| + \|v-a\| + \|u-v\|,
\]
so taking the infimum over $a$ on the right-hand side we obtain $\mu_d(u) \le \mu_d(v)+\|u-v\|$, and similarly, $\mu_d(v) \le \mu_d(u)+\|u-v\|$. Hence $|\mu_d(u) - \mu_d(v)| \le \|u-v\|$.

\begin{customthm}{2.7} \label{lem2.7}
For each $d>0$ there exists a locally Lipschitz continuous vector field $H: J^d\setminus K\to E$ with the following properties: \\
(i) $\|H(u)\|\le 1$.  \\
(ii) $\langle J'(u), H(u)\rangle  > \frac12 \mu_d(u)$. \\
(iii) $H$ has locally compact support, i.e.\ for each $u_0\in J^d\setminus K$ there 
exists a neighbourhood $U_0$ of $u_0$ in $J^d\setminus K$ and $R>0$ such that $u(x)=0$ for all $|x|\ge R$ and $u\in U_0$. \\
(iv) $H$ is odd in $u$.
\end{customthm}

\begin{proof}
Let $\wt u\in J^d\setminus K$. As $\|J'(\wt u)\| \ge \mu_d(\wt u)$, there exists $\wt v\in C_0^\infty(\rn)$, $\|\wt v\|=1$ such that $\langle J'(\wt u),\wt v\rangle >  \frac12\mu_d(\wt u)$ (this also holds if $\|J'(\wt u)\|=\infty$). Since $\wt v$ has compact support, $u\mapsto \langle J'(u), \wt v\rangle$ is continuous according to Lemma 2.2, and since also $\mu_d$ is  continuous, $\langle J'(u),\wt v\rangle > \frac12 \mu_d(u)$ for all $u$ in a neighbourhood $W(\wt u)$ of $\wt u$ in $J^d\setminus K$. Clearly, $(W(\wt u))$ is a covering of $J^d\setminus K$. Since $E$ is metric and hence paracompact, we can find a locally finite refinement $(W_j)$ and points $v_j\in C_0^\infty(\rn)$ such that $\langle J'(u),v_j\rangle > \frac12\mu_d(u)$ for $u\in W_j$. 
Let now
\[
H_1(u) := \sum_{j=1}^\infty\rho_j(u)v_j \quad \text{and} \quad H(u) := \frac12(H_1(u)-H_1(-u)),
\]
where $(\rho_j)$ is a Lipschitz continuous partition of unity subordinate to the cover $(W_j)$ of $J^d\setminus K$. It is easy to see that (i)-(iv) hold. 
\end{proof}

%\medskip
\noindent
The change in Lemma 2.7 prompts some (slight) changes in the proofs of Lemmas 2.13 and 2.14. First of all, the flow $\eta$ in (2.4) should be considered for $u\in J^d\setminus K$ instead of $u\in D(J)\setminus K$. Likewise, the flow $\wt\eta$ on p.\ 595 should be considered for $u\in J^{d+2\eps_0}\setminus K$.  On p.\ 592, line 17, $-z(\eta(t,u))$ should be replaced by $-\frac12\mu_d(\eta(t,u))$ and similarly, at the end of p.\ 592 in the definition of $\kappa_n$, $z(\eta(s,u))$ should read $\frac12\mu_d(\eta(s,u))$. Then $\kappa_n\to 0$ as before and by (ii) of Lemma \ref{lem2.7} and the properties of $\mu_d$, there exist $s_n^1\in[t_n,t_n^1]$ and $u_n^1$ such that $u_n^1-\eta(s_n^1,u)\to 0$, $J'(u_n^1)\to 0$ and $J(u_n^1)\le d$. Similarly we get $s_n^2\in[t_n^1,t_{n+1}]$ and $u_n^2$ such that $u_n^2-\eta(s_n^2,u) \to 0$, $J'(u_n^2)\to 0$ and $J(u_n^2)\le d$. Now the proof continues as in \cite{sqsz}.

A similar change must be made in the proof of Lemma 2.14. In the definition of $\tau$, $z(u)$ should be replaced by $\frac12\mu_{d+2\eps_0}(u)$. If $\tau=0$, we find a sequence $w_n^1\in J^{d+2\eps_0}_{d-2\eps_0}\cap U_\delta(K_d)\setminus U_{\delta/2}(K_d)$ with $\mu_{d+2\eps_0}(w_n^1) \to 0$ and then $u_n^1$ such that $u_n^1-w_n^1 \to 0$, $J'(u_n^1)\to 0$ and $J(u_n^1) \le d+2\eps_0$. The rest of the proof is the same as in \cite{sqsz}.

\medskip

At the end of the proof of Theorem 1.1 it is not clear that $T$ exists as claimed (the reason being the lack of continuity of $J$). On the other hand, $J$ is continuous on sets having compact support, see Lemma 2.2. In order to take advantage of this fact we need to slightly redefine the family $\mathcal{H}$ on p. 594.
Let 
\begin{align*}
\ch := \{h: D(J)\to E,\ \text{$h$ odd homeomorphism onto $h(D(J))$, $J(h(u))\le J(u)$} \\
\text{for all $u\in D(J)$ and if $A\in\Sigma$ has compact support, then so does $h(A)$}\}.
\end{align*}
Since for each $k\ge 1$ there exists  a $k$-dimensional subspace $E_k\subset C_0^\infty(\rn)$, according to Lemma 2.16 there are sets $A\in\Sigma$ having compact support and pseudoindex $i^*(A)\ge k$. Therefore the numbers
\[
d_k := \inf\{\sup_{u\in A} J(u): i^*(A)\ge k \text{ and $A$ has compact support}\}
\]
are well defined for all $k\ge 1$. Moreover, $d_k\ge b$ by Lemma 2.15.

Let $\wt\eta$ (with $u\in J^{d+2\eps_0}\setminus K$) be as in the proof of Theorem 1.1. We must show that there exists $T>0$ such that $J(\wt\eta(T,u)) < d-\eps$ for all $u\in B := A\setminus U$. It is clear by construction of  the vector field $H$ that $\wt\eta([0,t]\times B)$ has compact support for any $t>0$. Given $u_0\in B$, we have $J(\wt\eta(T_0,u_0)) < d-\eps$ for some $T_0>0$. Since the restriction of $J$ to $\wt\eta([0,T_0]\times B)$ is continuous according to Lemma 2.2, we can find a neighbourhood $V_0$ of $u_0$ in $B$ such that $J(\wt\eta(T_0,u)) < d-\eps$ for all $u\in V_0$. Now we use compactness in order to find a finite covering $V_i$ of $B$ and take $T := \max_i T_i$, where $T_i$ is such that $J(\wt\eta(T_i,u)) < d-\eps$ for all $u\in V_i$. 

A similar change needs to be made in the proof of Theorem 1.2: $\Gamma$ should consist of paths $\alpha$ which have compact support. See also Remark 3.8 in \cite{jisz}.

\end{document}